\documentclass[preprint,12pt]{elsarticle}

\usepackage{amssymb}
\usepackage{amsmath}

\usepackage[linesnumbered,lined,boxed,commentsnumbered]{algorithm2e}
 \usepackage{qtree}
 \usepackage{xcolor}
 \usepackage{newfloat}
\usepackage{tikz}
\usetikzlibrary{arrows.meta}
\DeclareFloatingEnvironment[fileext=lod]{diagram}

\RestyleAlgo{ruled} 

\newcommand{\bbZ}{\mathbb{Z}}
\newcommand{\bbN}{\mathbb{N}}

\newtheorem{theorem}{Theorem}

\newtheorem{definition}[theorem]{Definition}
\newtheorem{exam}[theorem]{Example}
\newtheorem{prop}[theorem]{Proposition}
\newproof{proof}{Proof}
\newtheorem{corollary}[theorem]{Corollary}

\begin{document}

\begin{frontmatter}



\title{Gapset Extensions, Theory and Computations} 

\author[label1]{Arman Ataei Kachouei} 
\author[label2]{Farhad Rahmati \corref{cor1}}
\affiliation[label1]{organization={Department of Mathematics and Computer Science, Amirkabir University of Technology},
            city={Tehran},
            country={Iran},
            mail={, arman.ataei@aut.ac.ir}}
\affiliation[label2]{organization={Department of Mathematics and Computer Science, Amirkabir University of Technology},
            city={Tehran},
            country={Iran},
            mail={, frahmati@aut.ac.ir}}
\cortext[cor1]{Corresponding author}

\begin{abstract}
In this paper we extend some set theoretic concepts of numerical semigroups for arbitrary sub-semigroups of natural numbers. Then we characterized gapsets which leads to a more efficient computational approach towards numerical semigroups and finally we introduce the extension of gapsets and prove that the sequence of the number of gapsets of size $g$ is non-decreasing as a weak version of Bras-Amorós's conjecture.

\end{abstract}



\begin{keyword}
Bras-Amorós's Conjecture  \sep Gapset \sep Gapset Extensions \sep  Numerical Semi Groups\sep Standard Basis 


\end{keyword}

\end{frontmatter}


\section{Introduction}
By a numerical semigroup we mean an additive submonoid of $\bbN_0$, the set of non-negative integers,  with finite complement in $\bbN_0$. This  finite complement is the gapset of the numerical semigroup and the size of the gapset is it's genus. The smallest element of $\bbN$ not in the gapset is named the multiplicity of  numerical semigroup. After this, we use the notation NS for numerical semigroup. A finitely generated submonoid of $\bbN_0$ is a NS if and only if it can be generated by a finite relatively prime set \cite[chp1]{Assi}. 
 \par
 Let $S$ be a NS with $m$ as it's multiplicity and $[i]_m^S$, for $0\leq i<m$, the set of elements of $S\setminus\{0\}$ that are in congruent with $i$ modulo $m$. The congruent modulo $m$ is an equivalence relation which induces $m$  partition sets on $\bbN_0$. In section \ref{section-1}, it has been proved that each NS, $S$, has a finite generating set $H=[m, h_{1}, \ldots, h_{m-1}]$ where $m = \min(S)$ and $h_{i}$ is the smallest element of $[i]_m^S$. This special unique generating set is called the standard basis for the given numerical semigroup \cite{Assi}. In this section we will extend set theoretic notions of numerical semigroups to arbitrary additive sub-semigroups. 
 \par
 So we can represent a NS by either one of it's set of generators or it's gapset, which are both finite subsets of $\bbN$, i.e. both generating sets and the gapsets can be used to determine uniquely a NS. We use a combination of these two approaches to study numerical semigroups. In section \ref{section-2}, we will characterize gapsets in a way that has more computational advantages with respect to other characterizations. In section \ref{section-3}, we will apply this characterization to compute the standard basis for a given NS. 
 \par
  In \cite{ELIAHOU2020105129} using congruent modulo $i$ partitioning, the authors describe gapsets and their special extensions that fix $m$, called $m$-extensions, which is exactly the same approach that we follow in this paper. In section \ref{section-4}, we will study gapset extensions and prove that the sequence $n_g$, the number of numerical semigroups with genus $g$, is a non-decreasing sequence. This proves \cite[Conjecture 2]{Kaplan}. This conjecture is a weak version of the Bras-Amorós's conjecture \cite{BrasAmoros}  that says $n_g\geq n_{g-1} + n_{g-2}$ for $g\geq 2$.

\section{Sub-Semigroups of $\bbN_0$}\label{section-1}
In this section, we extend some set theoretic notions of NS to arbitrary additive sub-semigroups of $\bbN_0$.
\begin{definition}
 Let $S$ be an additively closed subset of $\bbN_0$ and $S^*=S\setminus\{0\}$. We define:
\begin{enumerate}
  \item  The multiplicity of S, $m(S) :=\min(S^*)$,
  \item  The gapset of S,  $G(S):=\bbN_0\setminus S$,
  \item if $G(S)$ is a finite set then, Frobenius element of S, $F(S):= \max(G(S))$,
  \item The pseudo-Frobenius elements of S,  $PF(S):=\{x\in G(S) \mid x+ S^* \subseteq S\}$
  \item $PF_*(S):=\{x\in G(S) \mid x + S^* \subseteq G(S) \}$
\end{enumerate}
\end{definition}
\begin{exam}
Let k be a natural number. We have:
\begin{itemize}
  \item $PF(k\bbN_0) = \emptyset $
  \item $PF_*(k\bbN_0) = G(k\bbN_0)$
  \item $PF(10\bbN_0 + 12\bbN_0 + 14\bbN_0 + 16\bbN_0 + 18\bbN_0) = \{ 2,4,6,8 \}$
  \item $PF_*(10\bbN_0 + 12\bbN_0 + 14\bbN_0 + 16\bbN_0 + 18\bbN_0) = \{ 2k +1 \mid k \in \bbN_0 \} $
\end{itemize}
\end{exam}
We denote the set of elements of $X\subseteq \bbN_0$ which are congruent with $i$ modulo m  as 
$[i]_m ^{X}$ and by $[i]_m $ we mean $[i]_m ^{\bbN_0}$. The set $PF_*(S)$  has the following elementary properties.
\begin{prop}
Let $S$ be an additively closed subset of $\bbN_0$ and $x\in PF_*(S)$, then
$[x]_{m(S)} \subseteq PF_*(S)$.
\end{prop}
\begin{proof}
Every element $x\in PF_*(S)$ is a gap element. By definition \ref{stm:gap-seprator}, all elements of $[x]_{m(S)}$  less than $x$, are gap elements. Since $x + S^* \subseteq S^c$ and $km +S^* \subseteq S^*$, hence $x+km +S^* \in S^c$ for every integer $k>0$. 
\end{proof}
\begin{corollary}\label{cor:PF-star-and-gapsets}
If $PF_*(S)$ is a non-empty set, then $G(S)$ is not a finite set. 
\end{corollary}
The converse of corollary \ref{cor:PF-star-and-gapsets} is also true. Considering $S\cup \{0\} $, by \cite[proposition 1]{Assi}, $\gcd(S)=d\neq 1$. 
Let $x\in S$, obviously $d\mid m$, and $d\mid x$, hence $d\mid r_x$ where $r_x$ is the remainder of $x$ modulo $m$. So $[1]_m \subseteq G(S)$. For 
all unit element $u\in \bbZ_m$ there exists an element $k\in \bbN$ such that $ku \overset{m}{\equiv} 1$. We claim that there is no element of 
$[u]_m$ in $S$. Then we conclude that $S$ does not contain any element of a unit class from $\bbZ_m$, i.e. every unit classes of $\bbZ_m$ is a subsets of $G(S)$. We will show that sum of a unit and a zero-devisor is a unit in $\bbZ_m$. Let $u\in \bbZ_m$ be a unit and $z\in \bbZ_m$ be a zero-devisor. By the remainder theorem, $u = k_u m + r_u$ and $z = k_z m + r_z$. If $u+z$ is a zero-devisor then, $r_u + r_z = m$. So
$\gcd(m,u) = \gcd(m, r_u) = \gcd(m, -r_u)= \gcd(m,m-r_u) = \gcd(m, r_z) = \gcd(m, z)$. This is a contradiction. So we have the following proposition. 
\begin{prop}
 $PF_*(S)$ is a non-empty set if and only if $G(S)$  is an infinite set.
\end{prop}
If $G(S)$ is not a finite set, considering $S\cup \{0\}$, $\gcd(S)=d \neq 1$, hence $S\subseteq d\bbN_0$. Let $g\in G(S)\setminus d\bbN_0$ then, obviously, $d\nmid (g+s)$ for every $s\in S$. So the following proposition is proved.
\begin{prop}
$PF(S) \subseteq  d\bbN_0$, where $d = \gcd(S)$.
\end{prop}

It is well known that every sub monoid of $\bbN_0$ is a finitely generated monoid \cite[chp1, Prop 2]{Assi}. The following direct proof of this theorem suggest the general theme of visualizing a numerical semigroup like \cite{Delgado2020}. 
\begin{theorem}\label{finitely-generated}
Let $S$ be an additively closed subset of $\bbN_0$. There is a finite subset $\{s_1, \ldots, s_k \} \subseteq S$ such that 
$S = s_1 \bbN_0 + \cdots + s_k \bbN_0$.
\end{theorem}
\begin{proof}
  As a consequence of the well-ordering principle, $S^* = S \setminus \{0\}$
  has the least element $m$. The congruence modulo $m$ is an equivalence relation on $S$. 
  Every equivalence class of this relation, as a subset of $\bbN_0$, has the least element. Let $H_m[S]$ be the set of all such elements, then it is easy to see that $S = \sum_{h\in H_m} h\bbN_0$.
\end{proof}
In the proof of theorem \ref{finitely-generated}, we can use any element $k\in S$ such that $k>m(S)$, instead of it's multiplicity. As a similar way to \cite{Delgado2020}, in this paper we use trees for visualizing numerical semigroups. For a given additively closed subset $S\subseteq\bbN_0$, by the phrase “the $H$-set of $S$" we mean $H_m[S]$ where $m=\min(S)$. The following example shows that we can represent a semigroup $S\subseteq\bbN_0$, with the $H_k[S]$ set for every $k\geq m$, that introduced in the proof of theorem \ref{finitely-generated}.
\begin{exam}
Let $S=4\bbN_0  + 5\bbN_0+ 7\bbN_0$ be a numerical semigroup. The following tree visualizes this semigroup with respect to both $H_5[S]$ and $H_4[S]$, where the red nodes are the gap elements and the green ones are the elements of $H_k$ for each cases. This visualization shows that $H_k$ acts like a boundary between $S$ and $G(S)$. 
\begin{diagram}[ht]\label{Tree-representation}
\centering
\resizebox{0.9\textwidth}{!}{
\qtreecenterfalse
\Tree[.$H_4[S]$
        [.$\color{red}0$ 
            [.$\color{green}4+0$ 
                [.$2\cdot 4+0$ 
                    [.$\vdots$ ]
                ]            
            ]   
        ]
        [.$\color{red}1$  
            [.$\color{green}4+1$ 
                [.$2\cdot 4+1$ 
                    [.$\vdots$ ]
                ]            
            ]
        ]
        [.$\color{red}2$  
            [.$\color{red}4+2$ 
                [.$\color{green}2\cdot 4+2$ 
                    [.$\vdots$ ]
                ]            
            ]
        ]
        [.$\color{red}3$  
            [.$\color{green}4+3$ 
                [.$2\cdot 4+3$ 
                    [.$\vdots$ ]
                ]            
            ]
        ]
        ]
\hskip 0.3in
\Tree[.$H_5[S]$
        [.$\color{red}0$ 
            [.$\color{green}5+0$ 
                [.$2\cdot 5+0$ 
                    [.$\vdots$ ]
                ]            
            ]   
        ]
        [.$\color{red}1$  
            [.$\color{red}5+1$ 
                [.$\color{green}2\cdot 5+1$ 
                    [.$\vdots$ ]
                ]            
            ]
        ]
        [.$\color{red}2$  
            [.$\color{green}5+2$ 
                [.$2\cdot 5+2$ 
                    [.$\vdots$ ]
                ]            
            ]
        ]
        [.$\color{red}3$  
            [.$\color{green}5+3$ 
                [.$2\cdot 5+3$ 
                    [.$\vdots$ ]
                ]            
            ]
        ]
        [.$\color{green}5-1$ 
            [.$2\cdot 5-1$ 
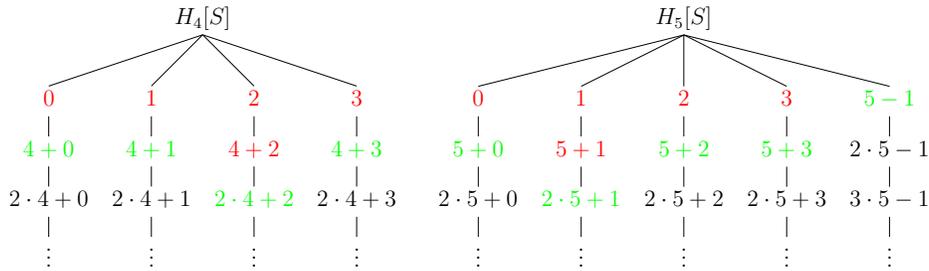
                [.$3\cdot 5-1$ 
                    [.$\vdots$ ]
                ]            
            ]
        ]
        ]
}
\vspace{4mm}
\caption{\centering Tree representations of $S=4\bbN_0  + 5\bbN_0+ 7\bbN_0$ with respect to $H_4$ and $H_5$}

\end{diagram}
\end{exam}
According to some authors, $H_m[S]$ is called the standard basis of $S$, whenever $S$ is a NS and  $m=\min(S)$ \cite[chapter 1]{Assi}.
\par
Let $S\subseteq\bbN_0$ be a monoid and $H$ be it's $H$-set. It is easy to verify the following simple statements. 
\begin{enumerate}\label{eq-pf}
  \item $G(S) = \overset{m-1}{\underset{i=1}{\bigcup}} \{ x\in [i]_m \mid x< \inf([i]_m^S) \}$ \footnote{We know that $\inf( \emptyset ) = \infty$.
  } \label{stm:gap-seprator}.
  \item $PF(S) \subseteq H-m$ 
  \item if $G(S)$ is a finite set then, $F(S) = \max(H-m)$. \footnote{Otherwise we can say that $S$ does not have a Frobenius number.}
\end{enumerate}
\begin{exam}
Let $S=a\bbN_0 + b\bbN_0 $ be a numerical semigroup where $\gcd(a,b) =1$ and $a<b$.
It is obvious that $b$ is a generator of $\bbZ_{a}$, hence $H=\{a, b, 2b, \ldots, (a-1)b\}$
is the standard basis for $S$ and, $(a-1)b -a = ab-(a+b)$ is the Frobenius number of S.
\end{exam}
\begin{exam}
Let $S = s_1 \bbN_0 + \ldots+ s_n\bbN_0\ $ be a NS, where, $s_1 < s_2 < \cdots < s_n$ and $\gcd(s_1, s_k) =1$ for some $k$. In this case, $s_1 \bbN_0 + s_k\bbN_0 + 2s_k\bbN_0 + \cdots + (s_1 - 1) s_k \bbN_0 \subseteq S$. So, for every $i>0$ we have $h_{i}\leq js_k$ for some $j<s_1$.
\end{exam}
In next sections we provide an algorithm that can compute the standard basis for every given numerical semigroup.

\section{Gapsets}\label{section-2}
A subset $G\subseteq\bbN_0$ is said to be a gapset, whenever it's a gapset of a numerical semigroup\footnote{One may extend this definition to an arbitrary sub-semigroup of $\bbN_0$. }. A gap element $x$ of a gapset $G$ is fundamental if $kx\notin G $ for all $k>1$\cite[chapter 3, section 5]{Rosales}. The set of fundamental gaps of $G$ is denoted by $FG(G)$.
Gapset characterization is a useful tool for dealing with computational aspects of numerical semigroups, for example we can find the following characterizations of gapsets in different resources.
\begin{itemize}
  \item\cite{ELIAHOU2020105129} A finite subset $G\subseteq \bbN_0$ is a gapset if and only if for all $x,y\in\bbN_0$ it contains $x$ or $y$ whenever $x+y\in G$.
  \item\cite[proposition 4.57]{Rosales} Let $X$ be a finite subset of $\bbN_0\setminus\{0\}$. The following conditions are equivalent.
       \begin{enumerate}
         \item $\{t\in\bbN_0\;:\; \exists x \in X \;\; t\mid x  \}$ is a gapset and $x\nmid y$ for all $x,y\in X$ such that $x\neq y$.
         \item There exists a gapset $G$ such that $FG(G) = X$.
       \end{enumerate}
       
\end{itemize}
These characterizations describe some aspects of the algebraic nature of gapsets, but they are not much convenient for computational purposes. Following simple, yet fundamental theorem gives another characterization of these structures which has some computational advantages over these characterizations.
\begin{theorem}\label{thm:gapset}
  Let $G$ be a non-empty proper subset of $\bbN_0$, $A= \bbN_0 \setminus G$, $m = \min(A)$, $h_i = \inf([i]_m^{A})$ for $i\in\{1,2,\ldots, m-1\}$ and $h_0 = m$. Then $A$ is a semigroup if and only if the following three conditions hold.
  \begin{enumerate}
    \item $\max([i]_m^{G}) < \min([i]_m^{A}$) for all $i\in \{ 1, \ldots, m-1 \} $ such that $h_i < \infty$ and $[m]_m^{G}$ is either empty or a singleton consists of zero.
    \item $h_i + h_j \geq h_{i+j} $ for all $i,j\in \{0,1, \ldots, m-1 \} $ such that $h_i, h_j <\infty$.
    \item If $h_i,h_j < \infty$ then $h_{i+j} < \infty$ for all $i,j \in \{0, \ldots, m-1 \}$. 
  \end{enumerate}
\end{theorem}
\begin{proof}
  Let $A$ be a semigroup. Based on theorem \ref{finitely-generated}, it is possible that for some  $i\in\{0, 1, \ldots, m-1\}$, $[i]_m^A = \emptyset$. Put $h_i= \infty$ for all such $i$. With these settings, all the three conditions are satisfied.\\ 
  For the converse, we need to show that $x+y \in A$ for every $x,y \in A$. Obviously two sets $[x]_m^A$ and $[y]_m^A$ are non-empty, hence $h_x, h_y < \infty$.  In $\bbZ_m$, we have the following equalities.
  \begin{equation}\label{eqn:modular-arithmetic}
    [ h_x + ‌h_y ]_m =[h_x]_m + [h_y]_m = [x]_m + [y]_m = [x+y]_m=[h_{x+y}]_m
  \end{equation} 
  $h_{x+y}<\infty$ follows from condition 3 and $x+y \geq h_{x+y}$ follows from the second condition, hence $x+y \in A$  by condition 1 and equation \ref{eqn:modular-arithmetic} .
\end{proof}
Theorem \ref{thm:gapset} describes every numerical sub-semigroup and numerical over-semigroup of an arbitrary numerical semigroup, i.e. it shows how to extend or reduce a gapset by adding or removing elements to another gapset.
\par
Let $G$ be a finite non-empty proper subset of $\bbN_0$ not containing zero. For $m\geq \min(\bbN\setminus G)$, we can compute one $H$-set for $G$ denoted by $H_m[G]=[m, h_1, \ldots, h_{m-1}]$ where $h_i=\min([i]_m^{\bbN_0\setminus G})$. For checking that $G$ is a gapset, based on theorem \ref{thm:gapset}, it is both necessary and enough to verify the following two conditions:
\begin{enumerate}
  \item $\max([i]_m^{G}) < \min([i]_m^{\bbN_0\setminus G}$) for all $i\in \{0, \ldots, m-1 \} $.
  \item $h_i+h_j \geq h_{i+j}$ for all two heads $h_i, h_j$ in $H_m[G]$\footnote{One may call them the "Hills" and the "Hill-set" instead of "Heads" and "Head-set".}.
\end{enumerate}
 
\par
For $0<k<m$, we have the set of relations, $R_{m,k} = \{ \{i,j\} \;:\; i+j \overset{m}{\equiv} k \}$. For an $H$-set $H_m[G]$ satisfying first condition of theorem \ref{thm:gapset} to be a gapset, it is enough to verify that $h_k \leq\min(\{h_i+h_j\; :\; \{i,j\}\in R_{m,k}\})$ for all $k$. This can reduce complexity of checking that if an $H$-set is a gapset or not. For example, for $m=5$ we can extract $R_{m,k}$ for all $k$ by computing the following table.
\begin{table}[ht]
\parbox{.3\linewidth}{
\centering
\begin{tabular}{|c|c|c|c|c|}
\hline
+ & 1 & 2                       & 3 & 4                        \\ \hline
1 & 2 & 3                       & 4 & 0                        \\ \hline
2 &   & 4                       & 0 & 1                        \\ \hline
3 &   &                         & 1 & 2                        \\ \hline
4 &   & {\color[HTML]{32CB00} } &   & {\color[HTML]{000000} 3} \\ \hline
\end{tabular}
\vspace{.2cm}
\caption{$m=5$}\label{table:5-relations}
}
\parbox{.3\linewidth}{
\centering
\vspace{.4cm}
\begin{tabular}{|c|c|c|c|}
\hline
+ & 1 & 2 & 3 \\ \hline
1 & 2 & 3 & 0 \\ \hline
2 &   & 0 & 1 \\ \hline
3 &   &   & 2 \\ \hline
\end{tabular}
\vspace{.2cm}
\caption{$m=4$}
}
\parbox{.3\linewidth}{
\centering
\vspace{.8cm}
\begin{tabular}{|c|c|c|}
\hline
+  & 1 & 2 \\ \hline
1 & 2 & 0 \\ \hline
2 &   & 1 \\ \hline
\end{tabular}
\vspace{.2cm}
\caption{$m=3$}
}
\end{table}
Table \ref{table:5-relations} shows that an $H$-set of multiplicity $5$ is a gapset if and only if all the following statements holds.
\begin{enumerate}
  \item $h_1 \leq \min(\{2h_3, h_4+h_2 \})$
  \item $h_2 \leq \min(\{ 2h_1, h_3+h_4 \})$
  \item $h_3 \leq \min(\{ h_1+h_2, 2h_4 \})$
  \item $h_4 \leq \min(\{ h_1+h_3, 2h_2 \})$
\end{enumerate}

\section{Computing the Standard Basis}\label{section-3}
Let $S= s_1\bbN_0+ \ldots+ s_k\bbN_0$ be a NS. In this section we will use theorem \ref{thm:gapset} to compute the standard basis for the given NS. Let $T = [m, h_1', \ldots, h_{m-1}' ]$ be a generating set for $S$ such that, $h_i '\in [i]_m^S$ and $m=\min(S)$. The following proposition shows that we can obtain the standard basis of $S$ from $T$ by updating $h_i'$s in a way that the second condition of theorem \ref{thm:gapset} holds.
\begin{prop}[correctness and finiteness]\label{prop:H-Mod}
Algorithm \ref{alg-HMode} computes the standard basis of $S$.\footnote{We used $\#$ for comment lines.}\\
\begin{algorithm}[ht]
\label{alg-HMode}
\SetKwInOut{Input}{input}\SetKwInOut{Output}{output}
\Input{$T = [m, h_1', \ldots, h_{m-1}' ]$,  A generating set for $S$ such that, $h_i '\in [i]_m^S$ and $m=\min(S)$.}
\BlankLine
\Begin{
\emph{\# \textit{isStdBasis}: a flag that controls second condition of theorem \ref{thm:gapset}}\;
$\textit{isStdBasis} \leftarrow \textbf{False}$\;
\While{$\textbf{not}\;\; \textit{isStdBasis}$}{
$\textit{isStdBasis} \leftarrow \textbf{True}$\;
\For{$i\leftarrow 1$ \KwTo $m-1$}{
\For{$j\leftarrow i$ \KwTo $m-1$}{
$index \leftarrow (i+j) \% m$\;
\If{$T[index] > T[i] + T[j]$}{
$T[index] \leftarrow T[i] + T[j]$ \;
$isStdBasis \leftarrow \textbf{False}$ \;
}}}
}
}
\caption{$T$-Modifier}\label{alg:H-Mod}
\end{algorithm}
\end{prop}
\begin{proof}
Let $T_k = [m, h_0^k, \ldots, h_{m-1}^k]$ be the result of this algorithm after $k$ iteration at the end of the $\textbf{While}$ loop. Based on theorem \ref{thm:gapset}, it is enough to show that this algorithm ends after finite iterations. Let assume that it does not end at all, hence according to the algorithm, there is at least an $i$ that the chain $h_i^{k_1} >h_i^{k_2} > h_i^{k_3} \cdots$ never stops. This contradicts the fact that $h_i'$ is finite for every $i$.
\end{proof}

The following proposition shows that we can obtain $T=[m, h_1', \ldots, h_{m-1}' ]$ from the generating set $\{s_1, \ldots, s_k\}$ of $S$.
\begin{prop}\label{prop:zm-generators}
Let $\mathbb{T}_0 = \{[s_1]_m, \ldots, [s_k]_m\}$ and $\mathbb{T}_l = \{ [s+t]_m \;:\; [s]_m, [t]_m \in\mathbb{T}_{l-1}  \}$ for $l\geq 1$. The following statements are equivalent.
\begin{enumerate}
  \item $\mathbb{T}_0$ is a generating set for $\bbZ_m$.
  \item There is an $l_0$ such that $\mathbb{T}_{l_0} = \bbZ_m$.
\end{enumerate}
\end{prop}
\begin{proof}
Assume that $\mathbb{T}_0$ is a generating set for $\bbZ_m$. $\mathbb{T}_l \subseteq \bbZ_m$ for every $l$ since $\mathbb{T}_0 \subseteq \bbZ_m$. Let $l_0$ be the first index that $\mathbb{T}_{l_0}$ is additively closed under the addition modulo $m$. Let $[j]_m\in \bbZ_m$. Since $\mathbb{T}_0$ is a generating set for $\bbZ_m$, we have $[j]_m = \sum_{i=1}^{k}c_i[ t_i]_m =\sum_{i=1}^{k}[c_i t_i]_m = [\sum_{i=1}^{k}c_i t_i]_m$ for some $c_i \in \bbN_0$. $\mathbb{T}_{l_0}$ is closed under addition modulo $m$, so $[j]_m \in \mathbb{T}_{l_0}$. The other side is obvious.
\end{proof}
As a result of proposition \ref{prop:H-Mod} and \ref{prop:zm-generators}, algorithm \ref{alg:SdtBasis} computes the standard basis of $S$.\\
\begin{algorithm}[H]
\SetKwInOut{Input}{input}\SetKwInOut{Output}{output}
\Input{$Gen=[s_1, \ldots, s_k ]$,  A generating set for the NS $S$.}
\Output{$T$, The standard basis of $S$.}
\BlankLine
\Begin{
    $m\leftarrow \min(s_1, \ldots, s_k)$\;
    \emph{\#initializing T as a list of m zeros (python list comprehension)}\\
    $T\leftarrow [0 \; \boldsymbol{for} \; i \; \boldsymbol{in} \; range(m) ]$ \\
    
    \For {$i\leftarrow 0$ \KwTo $m-1$}{  
    
            $T[i] \leftarrow \inf(\{ t\in Gen \;:\; t\overset{m}{\equiv} i \})$
    }
    \emph{\# \textit{isStdBasis}: a flag that controls second condition of theorem \ref{thm:gapset}}\\
    $\textit{isStdBasis} \leftarrow \textbf{False}$\\
    \While{$\textbf{not}\;\; \textit{isStdBasis}$}{
    $\textit{isStdBasis} \leftarrow \textbf{True}$\\
        \For{$i\leftarrow 1$ \KwTo $m-1$}{
            \For{$j\leftarrow i$ \KwTo $m-1$}{
                 
                $index \leftarrow (i+j) \% m$\\
                 \If{$T[index] > T[i] + T[j]$}{
                    $T[index] \leftarrow T[i] + T[j]$ \\
                    $isStdBasis \leftarrow \textbf{False}$ \\
                 }
            }
        }
    }
    \Return T
}
\caption{StdBasis}\label{alg:SdtBasis}
\end{algorithm}
\par
Let $a,b \in \bbZ$, according to \cite[chapter 1]{Assi} $a\leq_S b$ is defined as $b-a\in S$. Using this partially order, for a given numerical semigroup $S$, $PF(S)$ is the set of maximal gap elements of $S$ with respect to $\leq_S$ denoted by $Max_{\leq_S}(\bbN_0\setminus S)$\cite[chapter 1, proposition 7]{Assi}. Hence $PF(S) = Max_{\leq_S}(H-m)$ knowing $PF(S)\subseteq H$ where $H$ is the standard basis of $S$.
\par
 In the following example we are going to compute Frobenius and pseudo-Frobenius elements of a given NS using it's standard basis that is computed by algorithm \ref{alg:SdtBasis}.

\begin{exam}
Let $S = 6\bbN_0+17\bbN_0+38\bbN_0$. Following tables are computed by the algorithm \ref{alg:SdtBasis}, applied on $S$. Hence $H=[6, 55, 38, 51, 34, 17]$ is the standard basis of $S$. We can easily extract other invariants of $S$ like Frobenius Element $F(S)=\max(H)-m(S)=49$ and pesodo-Frobenius elements $PF(S)=Max_{\leq_S}(H-m)=\{49, 45\}$. 
\begin{table}[h]
\centering
\begin{tabular}{|c|c|c|}
\hline
           & $h_2 = 38$       & $h_5=17$          \\ \hline
$h_2 = 38$ & $76 \in [4]_6^S$ & $55\in [1]_6^S$   \\ \hline
$h_5=17$   &                  & $34 \in [4]_6^S $ \\ \hline
\end{tabular}
\caption{\centering Computed in first iteration\label{table-alg-out}}
\vspace{.5cm}
\begin{tabular}{|c|c|c|c|c|}
\hline
                                 & \multicolumn{1}{l|}{$h_1 = 55$} & $h_2 = 38$       & \multicolumn{1}{l|}{$h_4 = 34$} & $h_5=17$          \\ \hline
\multicolumn{1}{|l|}{$h_1 = 55$} & $110 \in [2]_6^S$               & $93\in [3]_6^S$  & $89\in [5]_6^S$                   & -                 \\ \hline
$h_2 = 38$                       & -                               & $76 \in [4]_6^S$ & -                               & $55\in [1]_6^S$   \\ \hline
\multicolumn{1}{|l|}{$h_4 = 34$} & -                               & -                & $68\in [2]_6^S$                 & $51\in [3]_6^S$   \\ \hline
$h_5=17$                         & -                               & -                & -                               & $34 \in [4]_6^S $ \\ \hline
\end{tabular}
\caption{\centering Computed in first and second iteration }
\vspace{.5cm}
\begin{tabular}{|c|c|c|c|c|c|}
\hline
           & $h_1 = 55$        & $h_2 = 38$       & $h_3=51$         & $h_4 = 34$      & $h_5=17$          \\ \hline
$h_1 = 55$ & $110 \in [2]_6^S$ & $93\in [3]_6^S$  & $106 \in [4]_6^S$ & $89\in [5]_6^S$ & -                 \\ \hline
$h_2 = 38$ & -                 & $76 \in [4]_6^S$ & $89\in[5]_6^S$   & -               & $55\in [1]_6^S$   \\ \hline
$h_3=51$   & -                  & -                 & -                & $85\in [1]_6^S$ & $68\in[2]_6^S$    \\ \hline
$h_4 = 34$ & -                 & -                & -                 & $68\in [2]_6^S$ & $51\in [3]_6^S$   \\ \hline
$h_5=17$   & -                 & -                & -                 & -               & $34 \in [4]_6^S $ \\ \hline
\end{tabular}
\caption{\centering Computed in third/last iteration}

\end{table}
\end{exam}
Based on theorem \ref{thm:gapset}, applying algorithm \ref{alg:SdtBasis} on a given set of generators for an additively closed subset $S\subseteq\bbN_0$, returns a list $T$ containing it's $H$-set. $H_m[S]$ consists of finite elements of $T$.
\section{Gapset Extensions}\label{section-4}
Let $G$ be a finite gapset of multiplicity $m$ and it's corresponding head list $H_m[G]=[m, h_1, \ldots, h_{m-1}]$.
Based on theorem \ref{thm:gapset}, all the elements of $[i]_m^{\bbN_0}$ which are less than $h_i$ belongs to $G$ for all $i$. Based on the terminalogy given in \cite{ELIAHOU2020105129}, let $G_i=[i, i+m, \ldots, h_i-m]$ be the list of elements of $[i]_m^{\bbN_0} \cap G$ in the increasing order, it is easy to see that $G = \cup_{i=0}^{m-1}G_i$ and $G_i \cap G_j =\emptyset$. Every $m-$extension $G^e$ of $G$ with head list $H_m^e[G^e] = [m, h_1^e, \ldots, h_{m-1}^e]$ is the result of adding elements of $[i]_m^{\bbN_0}\cap [h_i, h_i^e]$ to $G$ for all $i\in\{1,2,\ldots,m-1\}$, where $[h_i, h_i^e]=\{x\in\bbN_0 \; :\; h_i \leq x\leq h_i^e\}$, in a way that the second condition of theorem \ref{thm:gapset} holds.
We have the following sets.
\begin{itemize}
  \item $N_g$: The set of all gapsets of genus $g$
  \item $N_{m,g}$: The set of all gapsets of multiplicity $m$ and genus $g$ 
\end{itemize}
Finding some invariants for $N_g$ or $N_{m,g}$, can help understanding the structure of these sets. The following proposition suggests that the sum of all of heads in the $H$-set of a gapset $G\in N_{m,g}$ is always equal to a constant.
\begin{prop}
Let $H=[m, h_1, \ldots, h_{m-1}]$ be an $H$-set of a gapset $G\in N_{m,g}$, then we have $\sum_{i=1}^{m-1} h_i = gm + \frac{m(m-1)}{2}$.
\end{prop}
\begin{proof}
Let $h_i = k_i m + i$. We know that $\mid G_i \mid = k_i$ hence $\sum_{i=1}^{m-1} h_i m = gm$ and $\sum_{i=1}^{m-1} i = \frac{m(m-1)}{2}$.   
\end{proof}
Describing the relation between $N_{m,g-1}$ and $N_{m,g}$ can be used to reveal the secret of how a gapset extends in numerical semigroups. For example, the next proposition shows that, not only $N_{m,g} \neq \emptyset$, but also the set of all extensions of every element of $N_{m,g}$ is a non-empty set for all $g\geq m-1$.
\begin{prop}\label{prop:add-remove-element-from-gapset}
Every element of $N_{m,g}$ has at least an extension in $N_{m, g+1}$ and every element of $N_{m,g+1}$ is an extension of an element of $N_{m,g}$ for all $g\geq m-1$.
\end{prop}
\begin{proof}
Let $H=[m, h_1, \ldots, h_{m-1}]$ be the $H$-set of an element $G\in N_{m,g}$. If there is no $i$ that $G\cup \{h_i + m\}$ becomes a gapset, then by theorem \ref{thm:gapset}, sum of every two elements of the finite set $\{h_1, \ldots, h_{m-1}\}$ belongs to it, which contradicts the finiteness of this set. The other statement follows from theorem \ref{thm:gapset}, with a similar proof.
\end{proof}
We can visualize $N_{m,g}$ as a tree, where each level of the tree consists of those $H_m$ sets that have the same genus. The first node of this tree is $H_m^0 = [m, m+1, m+2, \ldots, 2m-1]$ (first level of the tree) which is the only element of $N_{m, m-1}$. By proposition \ref{prop:add-remove-element-from-gapset}, every node in the next level of the tree can be obtained from a node in the last level, by adding an element $(m+i)+m$ to $G_i$ for some $i\in\{1,2,\ldots, m-1\}$, in a way that the second condition of theorem \ref{thm:gapset} holds, in this case it is satisfied for all $i$. Hence in the second level of this tree, we have $n_{m,m}=m-1$ gapsets. Two vertices in two consecutive levels are connected if and only if one of them is extension of the other. There is no edge between nodes of every two elements in two non-consecutive levels or in the same level. By continuing this process, this tree is a visualization of $N_{m,g}$ for all $g$. In this tree, By theorem \ref{thm:gapset}, maximum number of nodes in the next level connected to a node $\Delta \leq m-1$ and, by proposition \ref{prop:add-remove-element-from-gapset}, minimum number of nodes in the next level connected to a node $\delta \geq 1$.
\begin{exam}
$n_{2,g} =1$ for $g\geq 1$, since $\delta = \Delta = 1$ in the tree of $N_{2,g}$.
\end{exam}
\begin{exam}\label{exam-N-3-g}
The $N_{3,g}$ tree up to the fifth level is as follows.
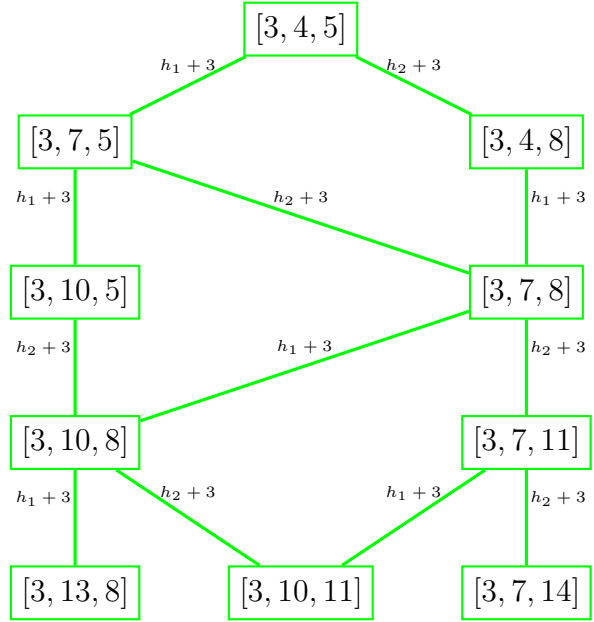
\begin{center}
\begin{tikzpicture}
\centering
\begin{scope}[every node/.style={thick,draw=green}]
    \node (n0) at (0,1.5) {$[3,4,5]$};
    \node (n1) at (3,0) {$[3,4,8]$};
    \node (n2) at (-3,0) {$[3,7,5]$};
    \node (n3) at (3,-2) {$[3,7,8]$};
    \node (n4) at (-3,-2) {$[3,10,5]$};
    \node (n5) at (3,-4) {$[3,7,11]$} ;
    \node (n6) at (-3,-4) {$[3,10,8]$} ;
    \node (n7) at (0,-6) {$[3,10,11]$} ;
    \node (n8) at (3,-6) {$[3,7,14]$} ;
    \node (n9) at (-3,-6) {$[3,13,8]$} ;
\end{scope}

\begin{scope}[>={[black]},
              every node/.style={},
              every edge/.style={above,draw=green,very thick}]
    \path [->] (n0) edge node {\tiny$h_2+3$} (n1);
    \path [->] (n0) edge node {\tiny$h_1+3$} (n2); 
    \path [->] (n1) edge node[xshift=12pt] {\tiny$h_1+3$} (n3);
    \path [->] (n2) edge node {\tiny$h_2+3$} (n3);
    \path [->] (n2) edge node[xshift=-12pt] {\tiny$h_1+3$} (n4);
    \path [->] (n3) edge node {\tiny$h_1+3$} (n6);
    \path [->] (n4) edge node[xshift=-12pt] {\tiny$h_2+3$} (n6);
    \path [->] (n6) edge node[xshift=-12pt] {\tiny$h_1+3$} (n9);
    \path [->] (n3) edge node[xshift=12pt] {\tiny$h_2+3$} (n5);
    \path [->] (n5) edge node[xshift=12pt] {\tiny$h_2+3$} (n8);
    \path [->] (n5) edge node[yshift=2pt] {\tiny$h_1+3$} (n7);
    \path [->] (n6) edge node[yshift=2pt] {\tiny$h_2+3$} (n7);
\end{scope}
\end{tikzpicture}
\end{center}
\end{exam}
\par
We know that, $n_{m,m-1}=1$ and $n_{m,m}=m-1$. 
The following theorem shows that if there are at least $m$ gapsets in $N_{m,g-1}$ then, $n_{m,g}\geq n_{m,g-1}$.
\begin{theorem}\label{thm:kaplans-inequivalences}
The following statements are True.
\begin{enumerate}
  \item If $n_{m,g-1} \geq m$ then, $n_{m,g} \geq n_{m,g-1}$.
  \item If $n_{m,g} \geq n_{m,g-1}$ then, $n_{m,g} \geq m-1$.
\end{enumerate}
\end{theorem}
\begin{proof}
Let $n_{m,g-1} \geq m$. Since $\Delta\leq m-1$, For every two different subsets $N_1, N_2 \subset N_{m,g-1}$ of size $m-1$, there are at least two different elements in $N_{m,g}$. So there are at least $n_{m,g-1} \choose m-1$ elements in $N_{m,g}$, hence ${n_{m,g-1} \choose m-1} \leq n_{m,g}$. If $n_{m,g} < n_{m,g-1}$ then, $n_{m,g-1} = m-1$ which is a contradiction.
\par
The second statement follows from the fact that, there are exactly $m-1$ extensions of $[m, m+1, m+2, \ldots,2m-1]$ which is the only element of $N_{m,m-1}$.
\end{proof}
Theorem \ref{thm:kaplans-inequivalences} shows that if we find a $g_0$ for which $n_{m,g_0}\geq m$ then, $n_{m,g} \geq n_{m,g-1}$ for every $g\geq g_0$. Based on theorem \ref{thm:kaplans-inequivalences}, Example \ref{exam-N-3-g} shows that $n_{3,g}\geq n_{3, g-1} \geq 2$ for $g\geq 3$. The following proposition proves that $g_0 = m+1
$ works fine for $m\geq 4$.
\begin{prop}
The set $N_{m,m+1}$ for $m\geq 4$ has at least $m$ elements.
\end{prop}
\begin{proof}
  We know that $n_{m,m}=m-1$ and the gapsets $H_m = [m, 2m+1, m+2, \ldots, 2m-1]$ and $H_m' = [m, m+1, m+2, \ldots, 3m-1]$ are elements of $N_{m, m}$. We prove that $H$ has at least two extensions and $H'$ has exactly $m-2$ extensions in $N_{m,m+1}$. By theorem \ref{thm:gapset}, $H'$ can be extends to $N_{m,m+1}$ from all of it's heads except possibly $h_{m-1}$, hence there are at least $m-2$ elements in $N_{m,m+1}$ and $H$ can be extends from $h_1$ and $h_2$ to an element of $N_{m, m+1}$. It is clear that none of these two extensions is equal to any of the $m-2$ extensions of $H'$, since the $h_{m-1}$ in these two extension is $2m-1$ while it is $3m-1$ for extensions of $H'$.
\end{proof}
\begin{corollary}\label{cor:none-decreasing-m-g}
$n_{m,g}$ is a non-decreasing sequence.
\end{corollary}
It is clear that $n_g = \sum_{m=2}^{g+1} n_{m,g}$, hence corollary \ref{cor:none-decreasing-m-g} proves that $n_g$ is a non-decreasing sequence.
  \bibliographystyle{elsarticle-num} 
  \bibliography{samplebib}



%
%
%
\end{document}